\documentclass[preprint,12pt]{elsarticle}

\usepackage{amsmath,amssymb,dsfont}
\usepackage{nomencl}
\makenomenclature

\usepackage{multirow}
\usepackage{graphicx}
\usepackage[percent]{overpic}
\usepackage{pict2e}

\usepackage[usenames,dvipsnames]{color}

\usepackage[bf,SL,BF]{subfigure}
 \makeatletter
 \newcommand{\setcaptype}[1]{%
 \renewcommand{\@captype}{#1}}
 \makeatother

\journal{Applied Mathematical Modelling}

\begin{document}

\begin{frontmatter}



\title{
A semi-analytical collocation method for solving multi-term variable-order time fractional partial differential equations}

\author[label1]{Xia Tian}
\author[label1]{S.Yu. Reutskiy\corref{cor1}\fnref{label4}}
\author[label1,label2,label3]{Zhuo-Jia Fu\fnref{label4}}
\address[label1]{Center for Numerical Simulation Software in Engineering \& Sciences, College of Mechanics and Materials, Hohai University, Nanjing, Jiangsu 210098, China}
\address[label2]{Institute of Continuum Mechanics, Leibniz University Hannover, Hannover 30167, Germany}
\address[label3]{State Key Laboratory of Mechanics and Control of Mechanical
	Structures, Nanjing University of Aeronautics \& Astronautics, Nanjing, Jiangsu 210016, China}
 \address[label4]{Corresponding authors: S.Yu. Reutskiy (sergiy.reutskiy2016@yandex.ru) and Z.J Fu (zhuojiafu@gmail.com)}


\begin{abstract}
This paper presents a novel semi-analytical collocation method to solve multi-term variable-order time fractional partial differential equations (VOTFPDEs). In the proposed method it employs the Fourier series expansion for spatial discretization, which transforms the original multi-term VOTFPDEs into a sequence of multi-term variable-order time fractional ordinary differential equations (VOTFODEs). Then these VOTFODEs can be solved by using the recent-developed backward substitution method. Several numerical examples verify the accuracy and efficiency of the proposed numerical approach in the solution of multi-term VOTFPDEs.
\end{abstract}


\end{frontmatter}


\section{Introduction}
In recent advances, instead of the standard partial differential equations, time fractional partial differential equations (TFPDEs) have been widely used to describe some anomalous natural processes in physics \cite{hilfer2000applications}, hydrology \cite{mandelbrot1968fractional,meerschaert2004finite} and finance \cite{scalas2000fractional}. According to various recent experimental results, it shows that constant-order TFPDEs even with multi-term time fractional derivatives cannot fully capture some time-dependent anomalous natural processes \cite{chechkin2005fractional,santamaria2006anomalous}. Therefore, variable-order TFPDEs have been introduced, in which the order of time fractional operator is a function of time. In this paper, we consider the following general form of multi-term variable-order time fractional partial differential equations (VOTFPDEs):%
\begin{multline}
D_{t}^{\alpha(t)}u(\mathbf{x},t)+\sum_{i=1}^{I}a_{i}\left( t\right)
D_{t}^{\alpha _{i}\left( t\right) }u(\mathbf{x},t)  \label{0.1} \\
=\left[ \sum_{i=I+1}^{m-1}a_{i}\left( t\right) D_{t}^{\alpha _{i}\left( t\right) }
\left(\sum_{i,j=1}^{d}\frac{\partial^{2} u(\boldsymbol{\mathrm{x}},t)}{\partial x_{i}\partial x_{j}}\right)\right]
 +f\left(
\mathbf{x},t\right) ,\text{ }\mathbf{x}\in \Omega \subset \mathbb{R}^{d},\text{ }0\leq t\leq T.
\end{multline}%
subjected to the following boundary conditions
\begin{equation}
u(\mathbf{x},t)=g_{1}\left(\mathbf{x},t\right),\text{ }\mathbf{x}\in \Gamma_1,\text{ }0\leq t\leq T.  \label{0.2}
\end{equation}%
\begin{equation}
\frac{\partial u(\mathbf{x},t)}{\partial \mathbf{n}}=g_{2}\left(\mathbf{x},t\right) ,\text{ }\mathbf{x}\in \Gamma_2,\text{ }0\leq t\leq T.  \label{0.3}
\end{equation}%
and initial conditions%
\begin{equation}
u(\mathbf{x},0)=h_{0}\left( \mathbf{x}\right) ,\text{ }\frac{\partial^{i} u(\mathbf{x},0)}{\partial t^{i}}%
=h_{i}\left( \mathbf{x}\right) ,\text{ }\mathbf{x}\in \Omega, i=1,2,...,m-1.  \label{0.4}
\end{equation}%
where $\mathbf{x}=(x_{1},x_{2},\cdots,x_{d})$, $\Omega=\prod_{i=1}^{d}{\left[0,L_{i} \right]}$, $\partial \Omega=\Gamma_1 \cup \Gamma_2$, $m-1< \alpha(t) \leq m$, $m-1-i\leq \alpha _{i}\left( t\right) \leq m-i,$ $i=1,...,m-1$ - are functions of $t$ in the interval $\left[ 0,T\right] $, $m \in \mathbb{N}^{+}$, $a_{i}\left( t\right) $, $g_{1}\left(\mathbf{x},t\right) $, $%
g_{2}\left(\mathbf{x},t\right) $, $%
f\left(\mathbf{x},t\right) $, $h_{0}\left( \mathbf{x}\right) $, $h_{i}\left( \mathbf{x}\right) $ are
known smooth enough functions. The definition of the variable-order time fractional derivative $D_{t}^{\alpha(t)}$ \cite{coimbra2003mechanics,zhao2015second} can be stated as follows
\begin{equation}
D_{t}^{\alpha(t)}F(t)=\left\{
\begin{array}{lll}
\frac{1}{\Gamma \left( m-\alpha(t) \right) }\int_{0}^{t}\frac{F^{\left( m\right)
}\left( \tau\right) }{\left( t-\tau\right) ^{\alpha(t) -m+1}}, &  & m-1<\alpha(t) <m, \\
F^{\left( m\right) }\left( t\right) , &  & \alpha(t) \equiv m,%
\end{array}%
\right.  \label{0.5}
\end{equation}%
where $\Gamma \left( z\right) $\ denotes the gamma function. According to this definition (\ref{0.5}), for the power functions we get:%
\begin{equation}
D_{t}^{\alpha(t)}t^{p}=\left\{
\begin{array}{lll}
0, \text{if }p\in \mathbb{N}\text{ and }p<m, \\
\frac{\Gamma \left( p+1\right) }{\Gamma \left( p+1-\alpha(t) \right) }t^{p-\alpha(t) },
 \text{if }p\in \mathbb{N}\text{ and }p\geq m\text{ or }p\notin
\mathbb{N}\text{ and }p>m-1,%
\end{array}%
\right.  \label{0.6}
\end{equation}%

It should be mentioned that Eq.(\ref{0.1}) is the general equation of various types of VOTFPDEs, such as the variable-order time-fractional diffusion equation \cite{zhuang2009numerical,sun2012finite}, the variable-order time-fractional telegraph equation \cite{hashemi2016numerical}, the variable-order time-fractional diffusion-wave equation \cite{zhang2012compact,liu2013numerical} as particular cases, to mention just a few.

Numerical simulation plays an important role on the investigation of the time fractional partial differential equations (VOTFPDEs). Nowadays, the finite difference methods (FDMs) are the popular and dominant numerical approaches for temporal and spatial discretization of the VOTFPDEs. Their convergence,accuracy,and stability have extensively been discussed in the literatures \cite{langlands2005accuracy,cui2009compact,lin2009stability}.

In addition, with traditional FDMs for temporal discretization, several highly-accurate and efficient numerical discretization methods have been introduced to spatial discretization of transformed time-independent partial differential equations, such as the Fourier series expansion method \cite{chen2010numerical}, the spectral method \cite{lin2007finite}, the finite element method \cite{li2011numerical}, the boundary element method \cite{katsikadelis2011bem}, and the radial basis function meshless collocation method \cite{brunner2010numerical,chen2010fractional,fu2013boundary}, to mention just a few. However, numerical investigation shows that the results are very sensitive to time discretization formulation in the FDM solution of the VOTFPDEs \cite{fu2015method}.

On the other hand, the M$\ddot{u}$ntz polynomials \cite{borwein1994muntz} of the form $\sum_{i=0}^{n}\beta_{i} z^{\lambda_{i}}$ with real coefficients $\beta_{i}$ can be constructed by the classical M$\ddot{u}$ntz-Sz$\acute{a}$sz Theorem, where $0\leq \lambda_{0}<\lambda_{1}< \cdots \to \infty$ and $\sum_{i=0}^{\infty}\lambda_{i}^{-1}=+\infty$. Based on the orthogonalization, the M$\ddot{u}$ntz-Legendre polynomials \cite{esmaeili2011numerical} can be derived, which has been widely used to solve different PDEs. Recently Mokhtary et al. \cite{mokhtary2016muntz} developed the M$\ddot{u}$ntz-Legendre tau method for fractional differential equations. One of the authors \cite{reutskiy2017new} applied the M$\ddot{u}$ntz polynomials to deal with the time discretization in the solution of constant-order TFPDEs. In this study we make a first attempt to use the M$\ddot{u}$ntz polynomials for the time discretization in the solution of multi-term variable-order TFPDEs.

This paper presents a semi-analytical collocation method based on the Fourier series expansion and the backward substitution method with the M$\ddot{u}$ntz polynomials to solve multi-term variable-order time fractional partial differential equations (VOTFPDEs). In the proposed numerical scheme we first transform the nonhomogeneous boundary conditions along the spatial coordinate into a homogeneous boundary condition. Then, we seek a solution of the nonhomogeneous problem by using the Fourier series expansion in conjunction with the radial basis function interpolation \cite{chen2014recent}. Due to the orthogonality of the trigonometric basis functions, we get a sequence of the multi-term variable-order time fractional ordinary differential equations (VOTFODEs) instead of the initial VOTFPDEs. We solve these equations by the use of the backward substitution method with the M$\ddot{u}$ntz polynomials. In general, the proposed scheme can be considered as a semi-analytical collocation method.

A brief outline of the paper is as follows. Section 2 describes the proposed semi-analytical collocation method for solving multi-term VOTFPDEs. In Section 3, the efficiency and accuracy of the proposed scheme are examined in comparison with the analytical solutions and some reference results under several benchmark examples. Finally, in Section 4, we shall draw some conclusions based on the numerical simulations.

\section{Semi-analytical collocation method for multi-term VOTFPDEs}\label{sec method}
In this section a semi-analytical collocation method based on the Fourier series expansion and the backward substitution method \cite{reutskiy2017new} with the M$\ddot{u}$ntz polynomials is introduced to solve multi-term VOTFPDEs.

To solve Eqs. (\ref{0.1})-(\ref{0.4}), we first transform the nonhomogeneous boundary conditions into homogeneous boundary conditions. Let us denote
\begin{equation}
u\left( \mathbf{x},t\right) = v\left(
\mathbf{x},t\right) +s\left( \mathbf{x},t\right) .  \label{2.1}
\end{equation}%
The function $s\left( \mathbf{x},t\right) $ can be determined by using different schemes introduced in Section 3, and let the function $v\left( \mathbf{x},t\right) $ satisfy the following equation%
\begin{multline}
D_{t}^{\alpha(t)}v(\mathbf{x},t)+\sum_{i=1}^{I}a_{i}\left( t\right)
D_{t}^{\alpha _{i}\left( t\right) }v(\mathbf{x},t)-\left[ \sum_{i=I+1}^{m-1}a_{i}\left( t\right) D_{t}^{\alpha _{i}\left( t\right) }
\left(\sum_{i,j=1}^{d}\frac{\partial^{2} v(\boldsymbol{\mathrm{x}},t)}{\partial x_{i}\partial x_{j}}\right)\right]  \label{2.2} \\
=f\left(
\mathbf{x},t\right)-D_{t}^{\alpha(t)}s(\mathbf{x},t)-\sum_{i=1}^{I}a_{i}\left( t\right)
D_{t}^{\alpha _{i}\left( t\right) }s(\mathbf{x},t)\\
-\left[ \sum_{i=I+1}^{m-1}a_{i}\left( t\right) D_{t}^{\alpha _{i}\left( t\right) }
\left(\sum_{i,j=1}^{d}\frac{\partial^{2} s(\boldsymbol{\mathrm{x}},t)}{\partial x_{i}\partial x_{j}}\right)\right]
=\Theta \left( \mathbf{x},t\right) \text{ }.
\end{multline}%
and the homogeneous boundary conditions
\begin{equation}
v(\mathbf{x},t)=0,\text{ }\mathbf{x}\in \Gamma_1,\text{ }0\leq t\leq T.  \label{2.3}
\end{equation}%
\begin{equation}
\frac{\partial v(\mathbf{x},t)}{\partial \mathbf{n}}=0 ,\text{ }\mathbf{x}\in \Gamma_2,\text{ }0\leq t\leq T.  \label{2.4}
\end{equation}%
The initial conditions (\ref{0.4}) are transformed into the following ones
\begin{equation}
v(\mathbf{x},0)=h_{0}\left( \mathbf{x}\right)-s(\mathbf{x},0) \equiv
v_{0}\left( \mathbf{x}\right)%
 ,\text{ }i=1,2,...,m-1.  \label{2.5}
\end{equation}%
\begin{equation}
\frac{\partial^{i} v(\mathbf{x},0)}{\partial t^{i}}%
=h_{i}\left( \mathbf{x}\right)-\frac{\partial^{i} s(\mathbf{x},0)}{\partial t^{i}} \equiv
v_{i}\left( \mathbf{x}\right)%
 ,\text{ }i=1,2,...,m-1.  \label{2.6}
\end{equation}%
We seek the solution of Eq.(\ref{2.2}) in the form of Fourier series%
\begin{equation}
v(\mathbf{x},t)=\sum_{n_{1},n_{2},\cdots,n_{d}=1}^{\infty }w_{n_{1},n_{2},\cdots,n_{d}}\left( t\right) \prod_{i=1}^{d} \sin \left( \frac{n_{i}\pi x_{i} }{L_{i}}\right).  \label{2.7}
\end{equation}%
As a result, we get a sequence of the multi-term VOTFODEs%
\begin{equation}
D^{\alpha(t)}w_{n_{1},n_{2},\cdots,n_{d}}\left( t\right)=\sum_{i=1}^{m-1}\beta _{i,n_{1},n_{2},\cdots,n_{d}}\left( t\right)
D^{\alpha_{i}(t)}w_{n_{1},n_{2},\cdots,n_{d}}\left( t\right)+\theta _{n_{1},n_{2},\cdots,n_{d}}\left( t\right) .\
\label{2.8}
\end{equation}%
with the following initial conditions:%
\begin{multline}
w_{n_{1},n_{2},\cdots,n_{d}}(0)=h_{0}^{n_{1},n_{2},\cdots,n_{d}}, \text{ }\label{2.9} \
\frac{\partial^{i}w_{n_{1},n_{2},\cdots,n_{d}}(0)}{\partial t^{i}}=h_{i}^{n_{1},n_{2},\cdots,n_{d}} \text{ }.
\end{multline}%
where $\beta _{i,n_{1},n_{2},\cdots,n_{d}}\left( t\right) =-a_{i}\left( t\right) ,$ $%
1\leq i\leq I$, $\beta _{i,n_{1},n_{2},\cdots,n_{d}}\left( t\right) =-\sum_{j=1}^{d} \frac{n_{j}^{2}\pi^{2} a_{i}\left( t\right)}{L_{j}^{2}} ,$ $I+1\leq i\leq m-1$, and
\begin{multline}
h_{0}^{n_{1},n_{2},\cdots,n_{d}}=\frac{2^{d}}{\prod_{i=1}^{d}L_{i}}\int_{0}^{L_{d}}\cdots \int_{0}^{L_{1}}v_{0}\left( \mathbf{x}\right)\prod_{i=1}^{d} \sin \left( \frac{n_{i}\pi x_{i} }{L_{i}}\right) dx_{1}\cdots dx_{d}, \text{ }\label{2.91} \\
h_{i}^{n_{1},n_{2},\cdots,n_{d}}=\frac{2^{d}}{\prod_{i=1}^{d}L_{i}}\int_{0}^{L_{d}}\cdots \int_{0}^{L_{1}}v_{i}\left( \mathbf{x}\right)\prod_{i=1}^{d} \sin \left( \frac{n_{i}\pi x_{i} }{L_{i}}\right) dx_{1}\cdots dx_{d}, \text{ }\\
\theta _{n_{1},n_{2},\cdots,n_{d}}\left( t\right) =\frac{2^{d}}{\prod_{i=1}^{d}L_{i}}\int_{0}^{L_{d}}\cdots \int_{0}^{L_{1}}\Theta \left( \mathbf{x},t\right)\prod_{i=1}^{d} \sin \left( \frac{n_{i}\pi x_{i} }{L_{i}}\right) dx_{1}\cdots dx_{d}.
\end{multline}%

It should be mentioned that Eqs.(\ref{2.8}) and (\ref{2.9}) are independent and
can be solved separately by using the backward substitution method (BSM) \cite{reutskiy2017new}.\ Let $\varphi _{k}\left( t\right) $ be some system of basis
functions on $\left[ 0,T\right] $ so that the right hand side of Eq.(\ref{2.8}%
)\ can be represented in the form of the following series

\begin{equation}
\sum_{i=1}^{m-1}\beta _{i,n_{1},n_{2},\cdots,n_{d}}\left( t\right)
D^{\alpha_{i}(t)}w_{n_{1},n_{2},\cdots,n_{d}}\left( t\right)+\theta _{n_{1},n_{2},\cdots,n_{d}}\left( t\right)=\sum_{k=1}^{\infty }q_{k}\varphi _{k}\left(
t\right).\label{2.11}
\end{equation}

Thus, the original Eq.(\ref{2.8}) can be rewritten in the form
\begin{equation}
D^{\alpha(t)}w_{n_{1},n_{2},\cdots,n_{d}}\left( t\right)=\sum_{k=1}^{\infty }q_{k}\varphi _{k}\left(
t\right) .  \label{2.12}
\end{equation}

Then we assume that for each $\varphi _{k}\left( t\right) $ there exists
$\phi _{k}(t)$ given in the explicit analytic form, which satisfies the
equation
\begin{equation}
D^{\alpha(t)}\phi _{k}(t)=\varphi _{k}\left( t\right) .
\label{2.13}
\end{equation}

In this paper we use the above-mentioned M$\ddot{u}$ntz polynomials
\begin{equation}
\phi _{k}\left( t\right) =t^{\delta _{k}},\text{ }\delta _{k}=\alpha+\delta
\left( k-1\right) \text{, }k=1,2,3,\cdots  \label{2.14}
\end{equation}%
\ as the basis functions in the BSM. The parameters are chosen as $0<\delta \leq 1$ and $\alpha=m $.

Using Eq.(\ref{0.6}), we gain:
\begin{equation}
\varphi _{k}(t)=\frac{\Gamma \left( \delta _{k}+1\right) }{\Gamma \left( \delta
_{k}+1-\alpha(t)\right) }t^{\delta _{k}-\alpha(t) }.  \label{2.15}
\end{equation}

Let us denote
\begin{equation}
\Phi _{k}\left( t\right) =\phi _{k}\left( t\right) +\sum_{i=0}^{m-1} c_{k,i}t^{i},
\label{2.16}
\end{equation}%
where the coefficients $c_{k,i},i=0,1,2,\cdots$ are determined in such a way that $%
\Phi _{k}\left( t\right) $ satisfies the homogeneous initial conditions (%
\ref{2.9}):%
\begin{equation}
\Phi _{k}\left( 0\right) =0,\text{ }\frac{\partial^{i}\Phi _{k}\left( 0\right) }{\partial t^{i}}=0\
\Rightarrow c_{k,i}=0,\text{ }i=0,1,2,\cdots  \label{2.17}
\end{equation}%

In the same way, let $\bar{w}_{n_{1},n_{2},\cdots,n_{d}}\left( t\right) $ be a solution of the
homogeneous equation%
\begin{equation}
D^{\alpha(t)}\bar{w}_{n_{1},n_{2},\cdots,n_{d}}\left( t\right)=0.  \label{2.18}
\end{equation}%
For $m-1<\alpha(t) \leq m$ any function of the form
\begin{equation}
\bar{w}_{n_{1},n_{2},\cdots,n_{d}}\left( t\right) =\sum_{i=0}^{m-1} c_{k,i}t^{i}  \label{2.19}
\end{equation}%
satisfies (\ref{2.18}). The coefficients $c_{k,i},i=0,1,2,\cdots$ are determined  in such
a way that $\bar{w}_{n_{1},n_{2},\cdots,n_{d}}$ satisfies the initial conditions (\ref{2.9}) of
the original VOTFODEs. In particular, when $m=2$, we have
\begin{multline}
\bar{w}_{n_{1},n_{2},\cdots,n_{d}}\left( 0\right) =h_{0}^{n_{1},n_{2},\cdots,n_{d}},\text{ }\frac{\partial \bar{w}_{n_{1},n_{2},\cdots,n_{d}}(0)}{\partial t}=h_{1}^{n_{1},n_{2},\cdots,n_{d}}\text{ }\label{2.20}\\
\Rightarrow c_{k,0}=h_{0}^{n_{1},n_{2},\cdots,n_{d}}\text{, }%
c_{k,1}=h_{1}^{n_{1},n_{2},\cdots,n_{d}}.
\end{multline}%

When $m=1$,\ we get only one parameter $c_{k,0}=h_{0}^{n_{1},n_{2},\cdots,n_{d}}$.

As a result, the sum
\begin{equation}
w_{n_{1},n_{2},\cdots,n_{d}}\left( t,\mathbf{q}\right) =w_{N}\left( t,\mathbf{q}\right) =\bar{w}_{N}\left( t\right) +\sum_{k=1}^{\infty
}q_{N,k}\Phi _{k}\left( t\right),\text{ }N=1,\cdots,\prod_{i=1}^{d}n_{i} \label{2.21}
\end{equation}%
satisfies the multi-term VOTFODE (\ref{2.8}) and the initial conditions (\ref{2.9}) with any choice of the parameters $q_{N,k}$. The approximate solution in the form of the truncated series
\begin{equation}
w_{N,K}\left( t,\mathbf{q}\right) =\bar{w}_{N}\left( t\right)
+\sum_{k=1}^{K}q_{N,k}\Phi _{k}\left( t\right)  \label{2.22}
\end{equation}%
satisfies the truncated equation%
\begin{equation}
D^{\alpha(t) }w_{N,K}(t,\mathbf{q})=\sum_{k=1}^{K}q_{N,k}\varphi
_{k}\left( t\right) .  \label{2.23}
\end{equation}

To get the unknowns $q_{N,1},...,q_{N,K}$ we substitute $w_{N,K}\left( t,\mathbf{q}%
\right) $ back into equation (\ref{2.11}):%
\begin{equation}
\sum_{k=1}^{K}q_{N,k}\left[ \varphi _{k}\left( t\right) -\sum_{i=1}^{m-1}\beta
_{i,N}\left( t\right) \Phi _{k}^{\alpha_{i}(t)}\left( t\right) %
\right] =\theta_{N} \left( t\right) +\sum_{i=1}^{m-1}\beta _{i,N}\left( t\right)
\bar{w}_{N}^{\alpha_{i}(t)}\left( t\right) ,  \label{2.24}
\end{equation}%
where we denote%
\begin{equation}
\Phi _{k}^{\alpha_{i}(t) }\left( t\right) \equiv D^{\alpha
_{i}(t) }\Phi _{k}\left( t\right) =\frac{\Gamma \left( \delta _{k}+1\right) t^{\delta _{k}-\alpha_{i}(t)}}{%
\Gamma \left( \delta _{k}+1-\alpha_{i}(t)\right) }+\sum_{j=0}^{m-1}c_{k,j}D^{\alpha_{i}(t) }t^{j}
\label{2.25},
\end{equation}%
\begin{equation}
\bar{w}_{N}^{\alpha_{i}(t)}\left( t\right) \equiv D^{\alpha_{i}(t) }\bar{w}_{N}\left( t\right),  \label{2.26}
\end{equation}%
in which we use Eq.(\ref{0.6}) to get $\Phi _{k}^{\alpha_{i}(t)  }\left(
t\right) $ and $\bar{w}_{N}^{\alpha_{i}(t)  }\left( t\right) $.

Applying the collocation procedure to equation (\ref{2.24}), we get the
system of linear algebraic equations:
\begin{equation}
\sum_{k=1}^{K}q_{N,k}\left[ \varphi _{k}\left( t_{j}\right)
-\sum_{i=1}^{m-1}\beta _{i,N}\left( t_{j}\right) \Phi _{k}^{\alpha_{i}(t) }\left( t_{j}\right) \right] =\theta_{N} \left( t_{j}\right)
+\sum_{i=1}^{m-1}\beta _{i,N}\left( t_{j}\right) \bar{w}_{N}^{\alpha_{i}(t)
}\left( t_{j}\right) ,\text{ }  \label{2.27}
\end{equation}%
where $t_{j}$, $j=1,...,N_{c}\ $are the Gauss-Chebyshev (GC) collocation
points%
\begin{equation}
t_{j}=\frac{T}{2}\left[ 1+\cos \left( \frac{\pi \left( 2j-1\right) }{2N_{c}}%
\right) \right]  \label{2.28}
\end{equation}%
located on $\left[ 0,T\right] $. We take the number of the collocation points $N_{c}$ twice as many as the number of free parameters $K$. As a
result, we obtain an overdetermined linear system which can be solved by the standard least squares procedure. Using (\ref{2.1}), (\ref{2.7}) and (\ref{2.22}), we get the final approximate solution $u_{N,K}\left( \mathbf{x},t\right) $.

\section{Numerical simulations}\label{sec num}

In this section, we present some numerical experiments to verify the performance of the semi-analytical collocation method described in the previous section. In the first, we verify the performance of the proposed backward substitution method (BSM) with M$\ddot{u}$ntz polynomials in the solution of the multi-term variable-order time fractional ordinary differential equations (VOTFODEs). Then numerical comparisons between the proposed results and the reference results \cite{sun2012finite,yang2012computationally,wei2012analysis} for 1D multi-term variable-order time fractional partial differential equations (VOTFPDEs) are made. Finally, the proposed method is applied to solve the VOTFPDEs under 2D computational domains.
In order to validate the numerical accuracy of the proposed scheme, the maximum absolute error $Merr$ at the final time instant $T$ and the relative error $Rerr$ in the range $[0,T]$ are defined as follow
\begin{equation}
Merr=\max_{1\leq i\leq N_{t}}\left\vert u_{\text{exact}}\left(
\mathbf{x}_{i},T\right) -u_{N,K}(\mathbf{x}_{i},T)\right\vert  \label{3.1}
\end{equation}%
\begin{equation}
Rerr=\frac{\sum_{i=1}^{N_{t}} {\sum_{j=1}^{K_{t}} {\left[u_{\text{exact}}\left(%
\mathbf{x}_{i},t_{j}\right) -u_{N,K}\left(\mathbf{x}_{i},t_{j}\right)\right]}^2}}{\sum_{i=1}^{N_{t}} {\sum_{j=1}^{K_{t}} {u_{\text{exact}}^2\left(%
\mathbf{x}_{i},t_{j}\right)}}}  \label{3.200}
\end{equation}
where $N_{t}$ the test points $\mathbf{x}_{i}$ are uniformly distributed inside $%
\Omega$ and $K_{t}$ the test time instants $t_{j}$ are uniformly distributed in the range $%
[0,T]$. To illustrate the efficiency of the proposed method, the following approximation order ($AO$) and convergence order ($CO$) are defined as follow
\begin{equation}
AO\left( K\right) =\frac{\log \left( Error\left( K\right)\right) }{\log \left( K^{-1}\right) }.\label{3.201}
\end{equation}%
\begin{equation}
CO\left( N\right) =\frac{\log \left( Error\left( N/2\right) /Error\left( N\right) \right) }{\log \left( 2\right) }=\log _{2}\left( \frac{%
Error\left( N/2\right) }{E\left( N\right) }\right) .\label{3.2}
\end{equation}%
where $Error$ can be $Merr$ or $Rerr$, and $CO\left( N\right) =p$ denotes that the error decreases in $2^{p}$ times when the number of the free parameters increases from $N/2$ to $N$. Thus, $CO$ can be used to characterize the efficiency of the proposed method. Unless otherwise specified, the free parameter $\delta=0.25$ is chosen in the following numerical implementation.

\textbf{Example 1}. Let us first consider the following multi-term variable-order time fractional ordinary differential equation (VOTFODE)%
\begin{equation}
D^{\alpha(t)}w\left( t\right)=\sum_{i=1}^{m-1}\beta _{i}\left( t\right)
D^{\alpha_{i}(t)}w\left( t\right)+\beta _{0}\left( t\right) w\left( t\right)+\theta \left( t\right) .\
\label{3.3}
\end{equation}%
To illustrate the efficiency of the backward substitution method (BSM) with the M$\ddot{u}$ntz polynomials for solving the multi-term VOTFODEs. Here $m=4$, $\alpha(t) =3.2+0.5 \sin t$, $\alpha_{1}(t)=0.1+0.5\sin t, \alpha_{2}(t)=1+\cos t, \alpha_{3}(t)=2+0.1 e^{t}$, $\beta_{0}(t)=1+t^2, \beta_{1}(t)=-\sin t, \beta_{2}(t)=-\cos t, \beta_{3}(t)=-e^{-t}$ and the right-hand function $\theta \left( t\right) $ corresponds to the following analytical solution $w_{%
\text{exact}}\left( t\right) =t^6+t^4+t^2+1$. And the following initial conditions can be also derived from the aforementioned analytical solution
\begin{equation}
w\left( 0\right)=1, \frac{\partial w(0)}{\partial t}=0, \frac{\partial^{2} w(0)}{\partial t^2}=2, \frac{\partial^{3} w(0)}{\partial t^3}=0 \text{ } .\
\label{3.4}
\end{equation}%
Tables \ref{tab1_1}-\ref{tab1_3} present the numerical relative errors $Rerr$ and approximation orders of the proposed BSM with various parameters $\delta$ and $K$ for solving Example 1 at $T=0.01, 1, 100$. Generally speaking, the numerical errors of Example 1 obtained by the proposed BSM decay with the increasing $\delta$ or $K$, and it always has high approximation orders with $AO>4$. It can be observed from Tables \ref{tab1_1}-\ref{tab1_3} that the BSM with few terms of M$\ddot{u}$ntz polynomials can achieve very accurate results ($Rerr<10^{-15}$), and then the numerical accuracy cannot improve with the increasing terms of M$\ddot{u}$ntz polynomials due to the effect of the machine epsilon $\epsilon=2.22 \times 10^{16}$ in the common-used double-precision floating point arithmetic.
\begin{table}[bh]
\caption{The numerical relative errors $Rerr$ and approximation orders of the proposed BSM in Example 1 at $T=0.01$.}
\begin{center}
\begin{tabular}{ccccccc}
\hline
& \multicolumn{2}{c}{$\delta=0.1$} & \multicolumn{2}{c}{$\delta=0.25$} & \multicolumn{2}{c}{$\delta=0.5$} \\
$K$ & $Rerr\left( K\right) $ & $AO\left( K\right) $ & $Rerr\left(
K\right) $ & $AO\left( K\right) $ & $Rerr\left( K\right) $ & $AO\left(
K\right) $ \\ \hline
$3$ & $4.05\text{E-}14$ & $28.071$ & $9.74\text{E-}15$ & $29.367$ & $%
1.57\text{E-}14$ & $28.933$ \\
$4$ & $6.63\text{E-}15$ & $23.550$ & $2.37\text{E-}15$ & $24.294$ & $%
3.53\text{E-}16$ & $25.666$ \\
$5$ & $9.00\text{E-}16$ & $21.526$ & $2.22\text{E-}16$ & $22.395$ & $%
1.28\text{E-}16$ & $22.737$ \\
$6$ & $6.99\text{E-}16$ & $19.477$ & $1.45\text{E-}16$ & $20.353$ & $%
1.28\text{E-}16$ & $20.423$ \\
$7$ & $2.37\text{E-}16$ & $18.489$ & $2.17\text{E-}16$ & $18.535$ & $%
1.28\text{E-}16$ & $18.805$ \\
$8$ & $1.45\text{E-}16$ & $17.537$ & $1.68\text{E-}16$ & $17.469$ & $%
1.28\text{E-}16$ & $17.598$ \\
$9$ & $1.94\text{E-}16$ & $16.466$ & $1.28\text{E-}16$ & $16.654$ & $%
1.28\text{E-}16$ & $16.654$ \\ \hline
\end{tabular}%
\end{center}
\label{tab1_1}
\end{table}

\begin{table}[bh]
\caption{The numerical relative errors $Rerr$ and approximation orders of the proposed BSM in Example 1 at $T=1$.}
\begin{center}
\begin{tabular}{ccccccc}
\hline
& \multicolumn{2}{c}{$\delta=0.1$} & \multicolumn{2}{c}{$\delta=0.25$} & \multicolumn{2}{c}{$\delta=0.5$} \\
$K$ & $Rerr\left( K\right) $ & $AO\left( K\right) $ & $Rerr\left(
K\right) $ & $AO\left( K\right) $ & $Rerr\left( K\right) $ & $AO\left(
K\right) $ \\ \hline
$3$ & $1.28\text{E-}02$ & $3.969$ & $2.14\text{E-}02$ & $3.497$ & $%
2.06\text{E-}02$ & $3.532$ \\
$4$ & $8.20\text{E-}03$ & $3.466$ & $1.80\text{E-}03$ & $4.550$ & $%
4.29\text{E-}04$ & $5.593$ \\
$5$ & $4.34\text{E-}04$ & $4.811$ & $1.16\text{E-}04$ & $5.629$ & $%
1.04\text{E-}16$ & $22.867$ \\
$6$ & $3.59\text{E-}04$ & $4.426$ & $2.14\text{E-}05$ & $6.000$ & $%
9.09\text{E-}17$ & $20.615$ \\
$7$ & $8.96\text{E-}05$ & $4.790$ & $2.93\text{E-}06$ & $6.548$ & $%
1.82\text{E-}16$ & $18.625$ \\
$8$ & $1.18\text{E-}05$ & $5.456$ & $1.20\text{E-}07$ & $7.663$ & $%
1.90\text{E-}16$ & $17.407$ \\
$9$ & $2.63\text{E-}06$ & $5.848$ & $1.45\text{E-}16$ & $16.599$ & $%
1.13\text{E-}16$ & $16.712$ \\ \hline
\end{tabular}%
\end{center}
\label{tab1_2}
\end{table}

\begin{table}[bh]
\caption{The numerical relative errors $Rerr$ and approximation orders of the proposed BSM in Example 1 at $T=100$.}
\begin{center}
\begin{tabular}{ccccccc}
\hline
& \multicolumn{2}{c}{$\delta=0.1$} & \multicolumn{2}{c}{$\delta=0.25$} & \multicolumn{2}{c}{$\delta=0.5$} \\
$K$ & $Rerr\left( K\right) $ & $AO\left( K\right) $ & $Rerr\left(
K\right) $ & $AO\left( K\right) $ & $Rerr\left( K\right) $ & $AO\left(
K\right) $ \\ \hline
$3$ & $6.10\text{E-}03$ & $4.636$ & $4.40\text{E-}03$ & $3.945$ & $%
2.20\text{E-}03$ & $5.560$ \\
$4$ & $1.30\text{E-}03$ & $4.771$ & $6.68\text{E-}04$ & $5.274$ & $%
1.25\text{E-}04$ & $6.481$ \\
$5$ & $2.32\text{E-}04$ & $5.200$ & $6.71\text{E-}05$ & $5.971$ & $%
4.71\text{E-}16$ & $21.928$ \\
$6$ & $3.73\text{E-}05$ & $5.691$ & $4.90\text{E-}06$ & $6.824$ & $%
4.68\text{E-}16$ & $19.700$ \\
$7$ & $5.56\text{E-}06$ & $6.218$ & $2.32\text{E-}07$ & $7.850$ & $%
3.12\text{E-}16$ & $18.347$ \\
$8$ & $1.26\text{E-}06$ & $6.532$ & $2.02\text{E-}08$ & $8.520$ & $%
5.26\text{E-}16$ & $16.919$ \\
$9$ & $5.08\text{E-}07$ & $6.596$ & $3.56\text{E-}16$ & $16.190$ & $%
1.47\text{E-}16$ & $16.591$ \\ \hline
\end{tabular}%
\end{center}
\label{tab1_3}
\end{table}

\textbf{Example 2}. Next consider the variable-order time fractional diffusion equation with homogeneous Dirichlet boundary conditions and zero initial condition%
\begin{equation}
D^{\alpha(t)}u\left(x, t\right)=\beta _{0}\left( t\right) \frac{\partial^{2} u\left(x, t\right)}{\partial x^2}+f \left(x, t\right), \text{ }0\leq x\leq L,\text{ }0\leq t\leq T,\
\label{3.5}
\end{equation}%
subjected to the following boundary conditions
\begin{equation}
u(0,t)=u(L,t)=0, \text{ }0\leq t\leq T,  \label{3.6}
\end{equation}%
and initial conditions%
\begin{equation}
u(x,0)=0 ,\text{ }0\leq x\leq L,  \label{3.7}
\end{equation}%
where $f \left(x, t\right)=\left(\frac{2}{\Gamma \left(3-\alpha(t)\right)}t^{2-\alpha(t)}+\frac{\beta _{0}\left( t\right) \pi^2 t^2}{L^2}\right)\sin\left(\frac{x \pi}{L}\right).$ The corresponding analytical solution can be stated as
\begin{equation}
u_{\text{exact}}(x,t)=t^2 \sin\left(\frac{x \pi}{L}\right).\label{3.8}
\end{equation}%
For easy comparison, the parameters are set as $\beta _{0}\left( t\right)=0.01$, $L=10$, $\alpha (t)=0.8+0.2 t/T$ in this example. This is a problem with single space harmonic. The approximate solution can be represented as $u(x,t)=w(t) \sin\left(\frac{x \pi}{L}\right)$. Then the very accurate results can be obtained by using the proposed semi-analytical collocation method with few terms of M$\ddot{u}$ntz polynomials ($K=5$). Table \ref{tab1_4} lists the maximum absolute errors $Merr$ of Example 2 at various time instants $T=0.1,0.2,0.3,0.4,0.5$ by using the proposed method in comparison with the Crank-Nicholson FDM \cite{sun2012finite} with $h=0.1$ and $\tau=0.01$. From Table \ref{tab1_4}, it can be found that the proposed method with $K=4$ can achieve the similar numerical accuracy to the Crank-Nicholson FDM \cite{sun2012finite} with $\tau=0.01$, and the errors does not increase with the time evolution. It reveals that the proposed method can avoid the error accumulation effect in the standard FDMs.

\begin{table}[bh]
\caption{The maximum absolute errors $Merr$ of the proposed method and the Crank-Nicholson FDM \cite{sun2012finite} with $h=0.1$ and $\tau=0.01$ at various time instants $T=0.1,0.2,0.3,0.4,0.5$ in Example 2.}
\begin{center}
\begin{tabular}{ccccc}
\hline
& \multicolumn{3}{c}{Proposed method ($\delta=0.25$)} & \multicolumn{1}{c}{Crank-Nicholson FDM \cite{sun2012finite}} \\
$T$ & $K=3$ & $K=4 $ & $K=5$ & $h=0.1$ and $\tau=0.01$ \\ \hline
$0.1$ & $2.03\text{E-}02$ & $7.43\text{E-}04$ & $%
3.47\text{E-}18$ & $4.18\text{E-}04$\\
$0.2$ & $9.20\text{E-}03$ & $8.17\text{E-}04$ & $%
1.39\text{E-}17$ & $6.98\text{E-}04$\\
$0.3$ & $1.33\text{E-}02$ & $1.58\text{E-}04$ & $%
2.78\text{E-}17$ & $5.68\text{E-}04$\\
$0.4$ & $1.17\text{E-}02$ & $2.70\text{E-}04$ & $%
5.55\text{E-}17$ & $2.33\text{E-}04$\\
$0.5$ & $7.60\text{E-}03$ & $4.19\text{E-}04$ & $%
5.55\text{E-}17$ & $2.03\text{E-}03$\\ \hline
\end{tabular}%
\end{center}
\label{tab1_4}
\end{table}

\textbf{Example 3}. Next consider the variable-order time fractional diffusion equation with homogeneous Dirichlet boundary conditions%
\begin{equation}
D^{\alpha(t)}u\left(x, t\right)=\frac{\partial^{2} u\left(x, t\right)}{\partial x^2}+f \left(x, t\right), \text{ }0\leq x\leq 1,\text{ }0\leq t\leq T,\
\label{3.9}
\end{equation}%
subjected to the following boundary conditions
\begin{equation}
u(0,t)=u(1,t)=0, \text{ }0\leq t\leq T,  \label{3.10}
\end{equation}%
and initial conditions%
\begin{equation}
u(x,0)=10 x^2 \left(1-x\right) ,\text{ }0\leq x\leq 1,  \label{3.11}
\end{equation}%
where $f \left(x, t\right)= 20 x^2 \left(1-x\right) \left(\frac{t^{2-\alpha(t)}}{\Gamma \left(3-\alpha(t)\right)}+\frac{t^{1-\alpha(t)}}{\Gamma \left(2-\alpha(t)\right)}\right)-20 {\left(t+1\right)}^2 \left(1-3 x\right).$ The corresponding analytical solution can be stated as
\begin{equation}
u_{\text{exact}}(x,t)=10 x^2 \left(1-x \right) {\left(t+1\right)}^2.\label{3.12}
\end{equation}%
For easy comparison, the parameters are set as $\alpha(t)=\frac{2+\sin(t)}{4}$, $T=1$ in this example. Table \ref{tab1_5} presents the maximum absolute errors $Merr$ of Example 3 at $T=1$ by using the proposed semi-analytical collocation method and the reference method \cite{yang2012computationally} with $h=0.005$. It can be observed from Table \ref{tab1_5} that under the same space discretization ($N=200$ or $h=0.005$) the proposed method with $K=5$ can achieve the similar numerical accuracy to the reference method \cite{yang2012computationally} with time stepping size $\tau=0.000625$.

\begin{table}[h]
\caption{The maximum absolute errors $Merr$ of the proposed method and the reference method in Example 3 at $T=1$.}
\label{tab1_5}
\begin{center}
\begin{tabular}{cccccccc}
\hline
\multicolumn{4}{c}{Proposed method($\delta=0.25$)} &  &  & \multicolumn{2}{c}{Reference method \cite{yang2012computationally}} \\
$K$ &$N=100$ & $N=200$ & $N=250$ &  &  & $\tau $
& $h=0.005$ \\ \hline
$4$ &$1.44\text{E-}04$ &$1.36\text{E-}04$ &$1.34\text{E-}04$ &  &  &  $0.01$
& $2.09\text{E-}04$ \\
$5$ &$6.25\text{E-}05$ &$8.07\text{E-}06$ &$2.05\text{E-}06$ &  &  &  $0.005$
& $8.54\text{E-}05$ \\
$6$ &$6.25\text{E-}05$ &$8.07\text{E-}06$ &$2.05\text{E-}06$ &  &  &  $0.0025$
& $3.49\text{E-}05$ \\
$7$ &$6.25\text{E-}05$ &$8.07\text{E-}06$ &$2.05\text{E-}06$ &  &  &  $0.00125$
& $1.44\text{E-}05$ \\
$8$ &$6.25\text{E-}05$ &$8.07\text{E-}06$ &$2.05\text{E-}06$ &  &  &  $0.000625$
& $5.98\text{E-}06$ \\ \hline
\end{tabular}%
\end{center}
\end{table}

\textbf{Example 4}. In this example, we consider the multi-term variable-order time fractional partial differential equation (VOTFPDE)%
\begin{equation}
D^{\alpha(t)}u\left(x, t\right)+\sum_{i=1}^{m-1}D^{\alpha_{i}(t)}u\left(x, t\right)=%
\frac{\partial^{2} u\left(x, t\right)}{\partial x^2}+f \left(x, t\right), \text{ }-1\leq x\leq 1,\text{ }0\leq t\leq T,\
\label{3.13}
\end{equation}%
subjected to the following boundary conditions
\begin{equation}
u(-1,t)=g_{1}(-1,t), u(1,t)=g_{1}(1,t)\text{ }0\leq t\leq T,  \label{3.14}
\end{equation}%
and initial conditions%
\begin{equation}
u(x,0)=\frac{\partial u\left(x, 0\right)}{\partial t}=0 ,\text{ }-1\leq x\leq 1,  \label{3.15}
\end{equation}%
where $f \left(x, t\right)$, $g_{1}(-1,t)$ and $g_{1}(1,t)$ correspond to the analytical solution $u_{\text{exact}}(x,t)=\left[1/cosh(x-0.1)+1/cosh(x+0.1)\right] t^2$. Here $s(x,t)$ can be easily derived as $s(x,t)=g_{1}(-1,t)+\frac{x}{2} \left[g_{1}(1,t)-g_{1}(-1,t)\right]$. We transform this problem into the range $[0,2]$ by substituting $y=x+1$, and then, employ the semi-analytical collocation method described in section 2. The parameters are set as $m=4$, $\alpha(t)=1.25+t^2/20$, $\alpha_{1}(t)=1.2+\cos (t)/20$, $\alpha_{2}(t)=1.15+t/20$, $\alpha_{3}(t)=1.1+\sin (t)/20$, $T=1$ in this example. Table \ref{tabl_6} displays the numerical relative error $Rerr(N)$ and convergence rate $CO$ of the proposed method with $K=4$. It can be found from Table \ref{tabl_6} that the proposed method has rapid convergence rate with $CO \geq 2.4$, and the numerical errors are insensitive to the free parameter $\delta$ in this example.

\begin{table}[h]
\caption{The numerical error $Rerr(N)$ and convergence rate $CO$ of the proposed method with $K=4$ in Example 4.}
\label{tabl_6}
\begin{center}
\begin{tabular}{ccccccc}
\hline
& \multicolumn{2}{c}{$\delta=0.1$} &
\multicolumn{2}{c}{$\delta=0.25$} &
\multicolumn{2}{c}{$\delta=0.5$} \\
$N$ & $Rerr(N)$ & $CO\left( N\right) $ & $Rerr(N)$ & $CO\left( N\right) $ & $Rerr(N)$ & $CO\left( N\right) $
\\ \hline
$10$ & $3.96\text{E-}05$ & - & $3.96\text{E-}05$ & $ - $ & $3.96\text{E-}05$ & $ - $
\\
$20$ & $6.14\text{E-}06$ & $2.691$ & $6.14\text{E-}06$ & $2.691$ & $6.14\text{E-}06$ & $2.691$
\\
$40$ & $1.12\text{E-}06$ & $2.457$ & $1.12\text{E-}06$ & $2.457$ & $1.12\text{E-}06$ & $2.457$
\\
$80$ & $1.69\text{E-}07$ & $2.725$ & $1.69\text{E-}07$ & $2.725$ & $1.69\text{E-}07$ & $2.725$
\\
$160$ & $2.32\text{E-}08$ & $2.869$ & $2.32\text{E-}08$ & $2.869$ & $2.32\text{E-}08$ & $2.869$
\\
$320$ & $2.99\text{E-}09$ & $2.955$ & $2.99\text{E-}09$ & $2.955$ & $2.99\text{E-}09$ & $2.955$
\\ \hline
\end{tabular}%
\end{center}
\end{table}

\begin{figure}
  \centering
  \includegraphics[width=0.8\textwidth]{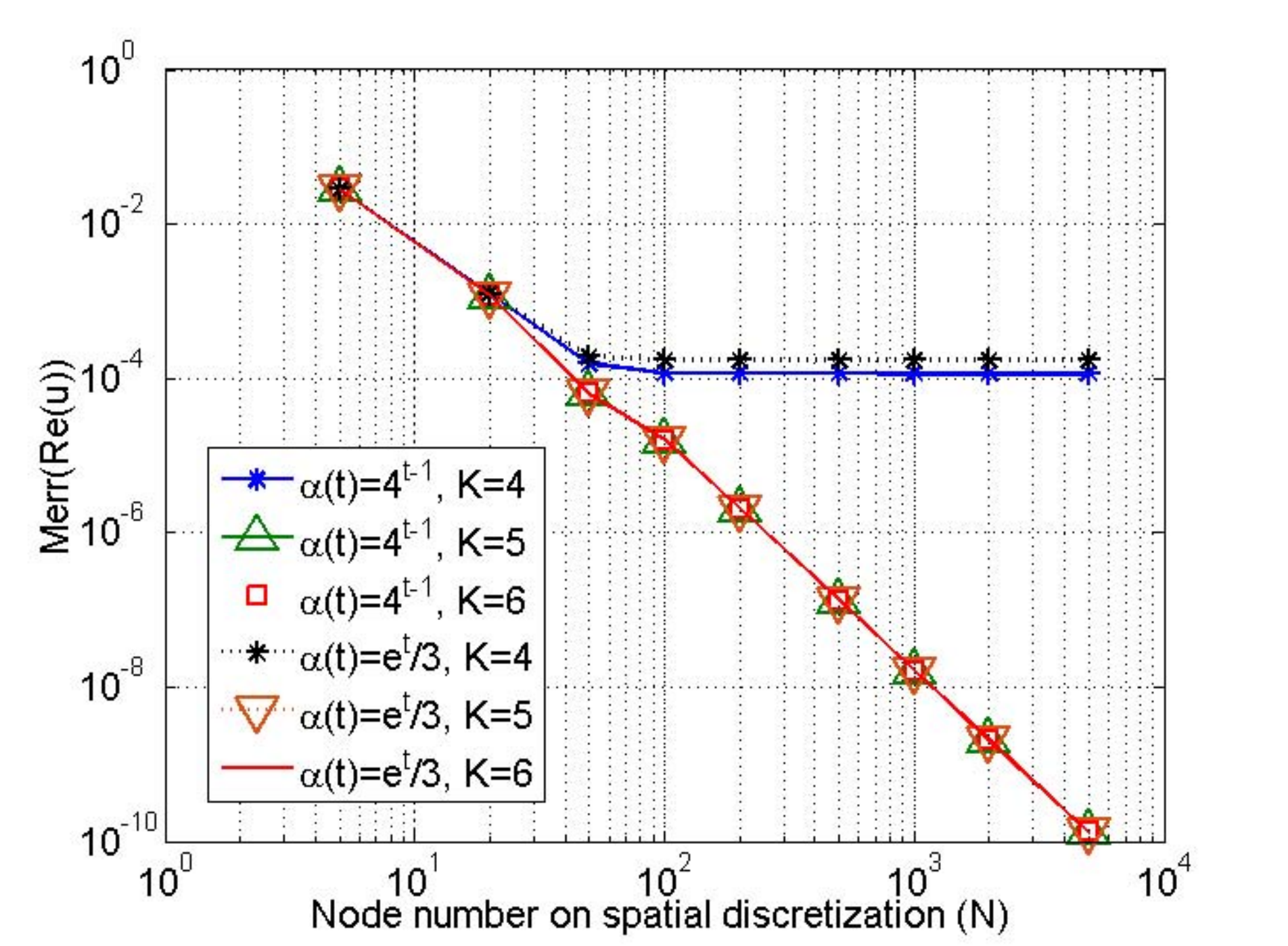}\\
  \caption{Numerical errors of the proposed method with respect to the node number on the spatial discretization $N$ in Example 5 with variable time fractional orders ($\alpha(t)=4^{t-1}$ and $e^{t}/3$) at $T=1$.}\label{Fig.1}
\end{figure}

\textbf{Example 5}. This example considers the variable-order time fractional Schr$\ddot{\text{o}}$dinger equation%
\begin{equation}
\sqrt{-1} D^{\alpha(t)}u\left(x, t\right)+\frac{\partial^{2} u\left(x, t\right)}{\partial x^2}=f \left(x, t\right), \text{ }0\leq x \leq 2 \pi,\text{ }0\leq t\leq T, \label{3.16}
\end{equation}%
subjected to the following boundary conditions
\begin{equation}
u(0,t)=u(2 \pi,t)=t^2, \text{ }0\leq t\leq T,  \label{3.17}
\end{equation}%
and initial conditions%
\begin{equation}
u(x,0)=0 ,\text{ }0\leq x\leq 2 \pi,  \label{3.18}
\end{equation}%
where $f \left(x, t\right)=-\frac{2 t^{2-\alpha(t)}}{\Gamma \left(3-\alpha(t)\right)} \sin(x)-t^2 \cos(x)+\sqrt{-1} \left(\frac{2 t^{2-\alpha(t)}}{\Gamma \left(3-\alpha(t)\right)} \cos(x)-t^2 \sin(x) \right)$. The corresponding analytical solution can be stated as
\begin{equation}
u_{\text{exact}}(x,t)=t^2 \left(\cos(x)+\sqrt{-1} \sin(x) \right).\label{3.19}
\end{equation}%
Here $s(x,t)$ can be easily derived as $s(x,t)=t^2$. First the time fractional order $\alpha(t)$ is set as a constant to make a numerical comparison between the proposed method and the implicit fully discrete local discontinuous Galerkin method in the literature \cite{wei2012analysis}. Table \ref{tabl_7} presents the maximum absolute errors $Merr$ of the proposed method and the reference method \cite{wei2012analysis} in Example 5 at T = 1 with different constant time fractional order $\alpha(t)=0.1,0.3,0.5$. Numerical comparison shows that under the equivalent spatial-temporal discretization the proposed method performs slight better than the reference method \cite{wei2012analysis}. Besides, it should be mentioned that the numerical error $Merr(Im(u))$ of the proposed method is less than $10^{-15}$ because that the Fourier series in Eq.(\ref{2.7}) includes the spatial part of the analytical solution. Then the proposed method is applied to Example 5 with variable time fractional orders ($\alpha(t)=4^{t-1}$ and $e^{t}/3$), the numerical errors are presented in Fig. \ref{Fig.1}. It can be found from Fig. \ref{Fig.1} that the proposed method only requires 5 terms of M$\ddot{u}$ntz polynomials ($K=5$) to obtain the enough accurate solutions in the time discretization.

\begin{table}[h]
\caption{The maximum absolute errors $Merr$ of the proposed method and the reference method \cite{wei2012analysis} in Example 5 at $T=1$.}
\label{tabl_7}
\begin{center}
\begin{tabular}{ccccccc}
\hline
& \multicolumn{3}{c}{Proposed method($\delta=0.25$)} & \multicolumn{3}{c}{Reference method \cite{wei2012analysis}} \\
$\alpha(t)$ &$(N,K)$ & $Re(u)$ & $Im(u)$ &$(N,K)$ & $Re(u)$ & $Im(u)$\\ \hline
$0.1$ &$(5,5)$ &$2.82\text{E-}02$ &$8.88\text{E-}16$ &$(5,5)$ &$3.18\text{E-}02$ &$3.12\text{E-}02$ \\
 &$(20,5)$ &$1.20\text{E-}03$ &$8.88\text{E-}16$ &$(10,10)$ &$4.11\text{E-}03$ &$3.89\text{E-}03$ \\
 &$(45,5)$ &$1.26\text{E-}04$ &$1.22\text{E-}15$ &$(15,15)$ &$1.22\text{E-}03$ &$1.22\text{E-}03$ \\
 &$(80,5)$ &$2.98\text{E-}05$ &$1.44\text{E-}15$ &$(20,20)$ &$5.19\text{E-}04$ &$5.16\text{E-}04$ \\\hline
$0.3$ &$(5,5)$ &$2.82\text{E-}02$ &$5.55\text{E-}16$ &$(5,5)$ &$3.19\text{E-}02$ &$3.11\text{E-}02$ \\
 &$(20,5)$ &$1.20\text{E-}03$ &$1.11\text{E-}15$ &$(10,10)$ &$4.07\text{E-}03$ &$3.85\text{E-}03$ \\
 &$(45,5)$ &$1.26\text{E-}04$ &$1.78\text{E-}15$ &$(15,15)$ &$1.18\text{E-}03$ &$1.18\text{E-}03$ \\
 &$(80,5)$ &$2.98\text{E-}05$ &$2.00\text{E-}15$ &$(20,20)$ &$4.83\text{E-}04$ &$4.84\text{E-}04$ \\ \hline
$0.5$ &$(5,5)$ &$2.82\text{E-}02$ &$5.55\text{E-}16$ &$(5,5)$ &$3.20\text{E-}02$ &$3.12\text{E-}02$ \\
 &$(20,5)$ &$1.20\text{E-}03$ &$1.11\text{E-}15$ &$(10,10)$ &$4.12\text{E-}03$ &$3.90\text{E-}03$ \\
 &$(45,5)$ &$1.26\text{E-}04$ &$1.33\text{E-}15$ &$(15,15)$ &$1.23\text{E-}03$ &$1.22\text{E-}03$ \\
 &$(80,5)$ &$2.98\text{E-}05$ &$1.78\text{E-}15$ &$(20,20)$ &$5.24\text{E-}04$ &$5.20\text{E-}04$ \\
 \hline
\end{tabular}%
\end{center}
\end{table}

\textbf{Example 6}. Next consider the following variable-order time fractional PDE%
\begin{equation}
D^{\alpha(t)}u\left(x, t\right)+D^{\alpha_{1}(t)}u\left(x, t\right)=\frac{\partial^2 u\left(x, t\right)}{\partial x^2}+f \left(x, t\right), \text{ }0\leq x\leq 1,\text{ }0\leq t\leq T, \label{3.20}
\end{equation}%
subjected to the following boundary conditions
\begin{equation}
u(0,t)=g_{1}(0,t),u(1,t)=g_{1}(1,t),\text{ }0\leq t\leq T,  \label{3.21}
\end{equation}%
and initial conditions%
\begin{equation}
u(x,0)=\frac{\partial u\left(x, 0\right)}{\partial t}=0 ,\text{ }0\leq x\leq 1,  \label{3.22}
\end{equation}%
where $f \left(x, t\right)$, $g_{1}(0,t)$, $g_{1}(1,t)$ can be easily derived by the analytical solution $u_{\text{exact}}(x,t)=t^2 e^{-100 (x-0.2)^2}$. Table \ref{tabl_8} displays the relative errors $Rerr(N)$ of the proposed method with $K=5$ in Example 6 with different variable time fractional orders in the range $[0,1]$. From Table \ref{tabl_8}, it can be found that the numerical results are in good agreement with the analytical solutions by using only 16 basis of Fourier series, and the proposed method have rapid convergence rate ($CO(N)>1.8$) with increasing basis of Fourier series.

\begin{table}[h]
\caption{The relative errors $Rerr$ of the proposed method in Example 6 with different variable time fractional orders in the range $[0,1]$.}
\label{tabl_8}
\begin{center}
\begin{tabular}{ccccccc}
\hline
& \multicolumn{2}{c}{$\alpha(t)=1.9+t/20$} & & & \multicolumn{2}{c}{$\alpha(t)=1.9+t/20$} \\
& \multicolumn{2}{c}{$\alpha_{1}(t)=1.6+\sin(t)/5$} & & & \multicolumn{2}{c}{$\alpha_{1}(t)=0.6+\cos(t)/5$} \\
$N$ &$Rerr(N)$ & $CO(N)$ &  & &$Rerr(N)$ & $CO(N)$ \\ \hline
$16$ &$7.09\text{E-}04$ & $-$ &  & &$7.09\text{E-}04$ & $-$ \\
$32$ &$1.19\text{E-}04$ & $2.578$ &  & &$1.19\text{E-}04$ & $2.578$ \\
$64$ &$3.31\text{E-}05$ & $1.845$ &  & &$3.31\text{E-}05$ & $1.845$ \\
$128$ &$2.11\text{E-}06$ & $3.967$ &  & &$2.11\text{E-}06$ & $3.967$ \\
$256$ &$5.04\text{E-}07$ & $2.068$ &  & &$5.04\text{E-}07$ & $2.068$ \\ \hline
\end{tabular}%
\end{center}
\end{table}

\textbf{Example 7}. Let us consider the following two-term variable-order time fractional wave-diffusion equation with damping in 2D square domain $\Omega=[0,1]^2$%
\begin{equation}
D^{\alpha(t)}u\left(\mathbf{x}, t\right)+\frac{\partial u\left(\mathbf{x}, t\right)}{\partial t}=%
\Delta {u\left(\mathbf{x}, t\right)}+f \left(\mathbf{x}, t\right), \text{ }\mathbf{x}=(x_{1},x_{2})\in \Omega,\text{ }0\leq t\leq T, \label{3.23}
\end{equation}%
subjected to the fully Dirichlet boundary conditions
\begin{equation}
u(\mathbf{x},t)=g_{1}(\mathbf{x},t),\text{ }\mathbf{x}\in \Gamma_1,\text{ }0\leq t\leq T,  \label{3.24}
\end{equation}%
and initial conditions%
\begin{equation}
u(\mathbf{x},0)=\frac{\partial u\left(\mathbf{x}, 0\right)}{\partial t}=0 ,\text{ }\mathbf{x}\in \Omega,  \label{3.25}
\end{equation}%
where $f \left(\mathbf{x}, t\right)=\left(\frac{6 t^{3-\alpha(t)}}{\Gamma \left(4-\alpha(t)\right)}+3 t^2-2 t^3 \right) e^{x_{1}+x_{2}}$ corresponds to the analytical solution $u_{\text{exact}}(\mathbf{x},t)=t^3 e^{x_{1}+x_{2}}$. Here $s(\mathbf{x}, t)$ can be represented by $s(\mathbf{x}, t)=\sum_{i=1}^{N_{b}}{\gamma_{i}\psi(\lVert \mathbf{x}-\mathbf{s}_{i}\rVert)}$, where $\psi$ represents the radial basis functions (RBFs),$\mathbf{s}_{i}$ is the RBF centers, and the unknown coefficients $\gamma_{i}$ can be determined by satisfying the boundary conditions \ref{3.24}. In the present implementation, $\psi$ is selected as the famous multi-quadric (MQ) RBFs $\psi(r)=\sqrt{r^2+c_{MQ}^2}$ with $c_{MQ}=4$, and the RBF centers $\{\mathbf{s}_{i}\} \in \Gamma$ with $N_{b}=4 \sqrt{N}-4$. Table \ref{tabl_9} presents the relative errors $Rerr$ of the proposed method with $K=4$ and $5$ in Example 7 with $\alpha(t)=1.85+\sin(t)/20$ in the range $[0,1]$. Numerical results show that the proposed method with $K=4$ and $N=25$ can perform very accurate results.

\begin{table}[h]
\caption{The relative errors $Rerr(u)$ and $Rerr(\frac{\partial u}{\partial x_{1}})$ of the proposed method in Example 7 with $\alpha(t)=1.85+\sin(t)/20$ in the range $[0,1]$.}
\label{tabl_9}
\begin{center}
\begin{tabular}{ccccc}
\hline
& \multicolumn{2}{c}{$K=4$} & \multicolumn{2}{c}{$K=5$} \\
$N$ &$Rerr(u)$ & $Rerr(\frac{\partial u}{\partial x_{1}})$ &$Rerr(u)$ & $Rerr(\frac{\partial u}{\partial x_{1}})$ \\ \hline
$25$ &$2.03\text{E-}05$ &$1.80\text{E-}03$ &$2.05\text{E-}05$ &$1.80\text{E-}03$ \\
$100$ &$7.73\text{E-}06$ &$6.62\text{E-}04$ &$7.58\text{E-}06$ &$6.62\text{E-}04$ \\
$225$ &$7.72\text{E-}06$ &$3.97\text{E-}04$ &$7.54\text{E-}06$ &$3.97\text{E-}04$ \\ \hline
\end{tabular}%
\end{center}
\end{table}

\begin{figure}
\setlength{\abovecaptionskip}{0pt}
\centering
\subfigure[Numerical solutions]{
\label{Fig.sub.(a)}
\includegraphics[width=0.45\textwidth]{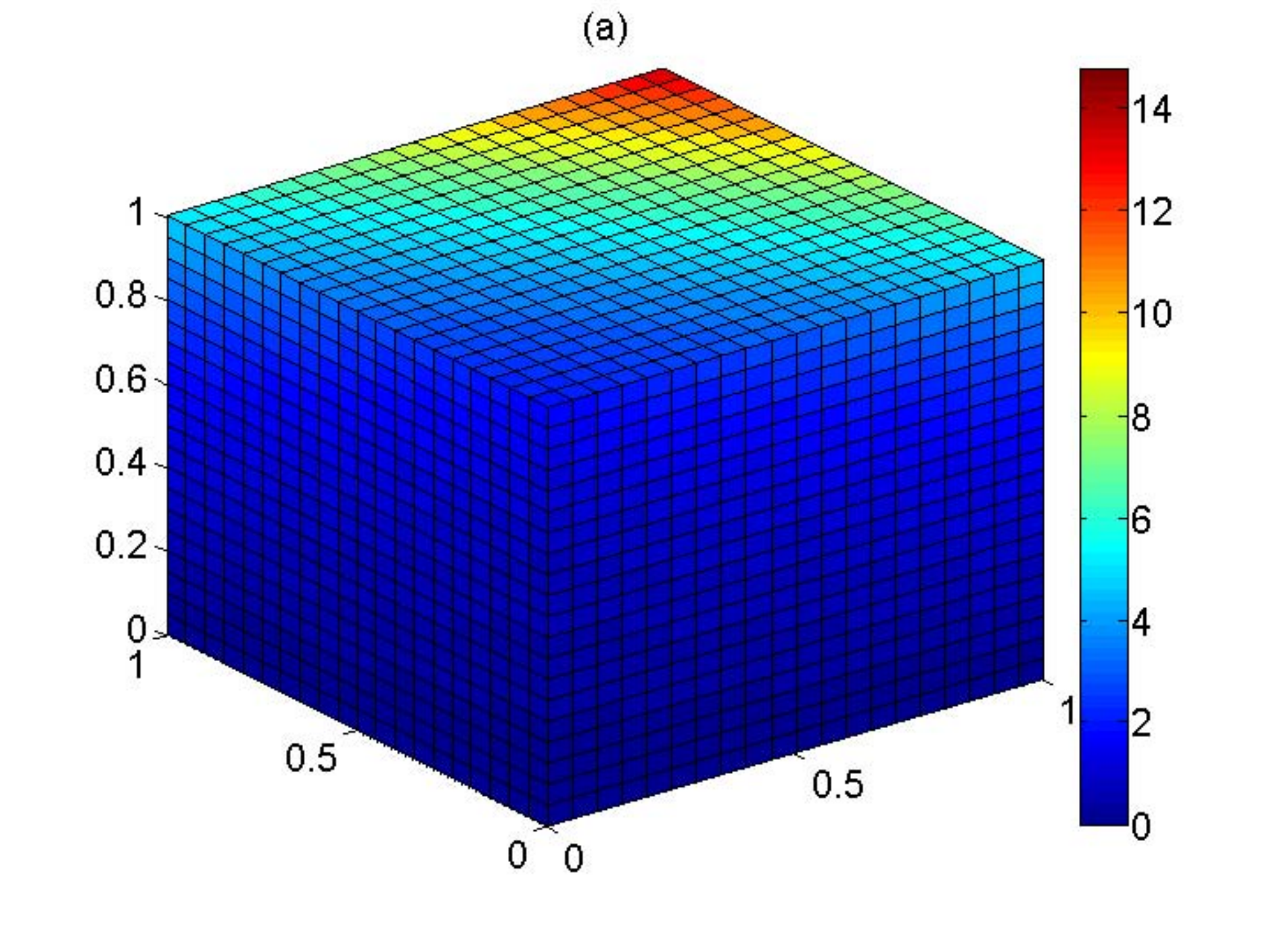}}
\subfigure[Numerical errors]{
\label{Fig.sub.(b)}
\includegraphics[width=0.45\textwidth]{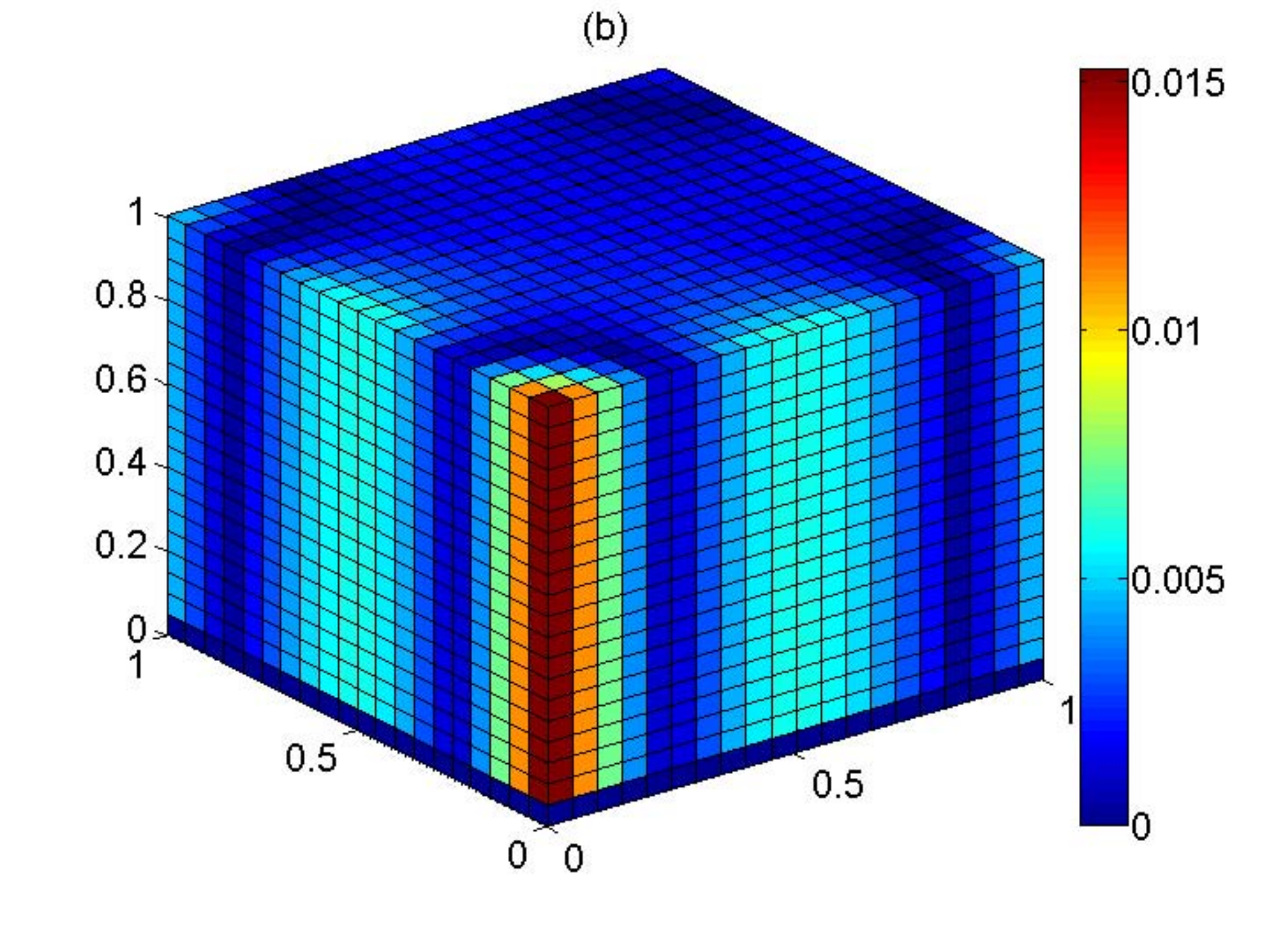}}\\
\caption{Numerical results of the proposed method with $K=5$ and $N=36$ in Example 8 with $\alpha=1.5$ and $\alpha(t)=1.4+t/10$ in the range $[0,1]$. }
\label{Fig.2}
\end{figure}

\textbf{Example 8}. Next consider the following variable-order time fractional wave-diffusion equation in 2D square domain $\Omega=[0,1]^2$%
\begin{equation}
D^{\alpha(t)}u\left(\mathbf{x}, t\right)+
\Delta^{2} {u\left(\mathbf{x}, t\right)}
=f \left(\mathbf{x}, t\right), \text{ }\mathbf{x}=(x_{1},x_{2})\in \Omega,\text{ }0\leq t\leq T, \label{3.26}
\end{equation}%
subjected to the following boundary conditions
\begin{equation}
u(\mathbf{x},t)=g_{1}(\mathbf{x},t), \Delta {u\left(\mathbf{x}, t\right)}=g_{2}(\mathbf{x},t),\text{ }\mathbf{x} \in \Gamma, \text{ }0\leq t\leq T,  \label{3.27}
\end{equation}%
and initial conditions%
\begin{equation}
u(\mathbf{x},0)=0, \frac{\partial u\left(\mathbf{x}, 0\right)}{\partial t}=e^{x_{1}+x_{2}} ,\text{ }\mathbf{x} \in \Omega,  \label{3.28}
\end{equation}%
where $f \left(\mathbf{x}, t\right)=e^{x_{1}+x_{2}}\left(\frac{\Gamma(4+\alpha) t^{3+\alpha-\alpha(t)}}{\Gamma \left(4+\alpha-\alpha(t)\right)}+ 4 t^{\alpha+3}+4 t \right)$, and $g_{1}(\mathbf{x},t)$, $g_{2}(\mathbf{x},t)$ can be easily derived from the analytical solution $u_{\text{exact}}(\mathbf{x},t)=e^{x_{1}+x_{2}} \left(t^{\alpha+3}+t\right)$.  Here $s(\mathbf{x}, t)$ can be also represented by $s(\mathbf{x}, t)=\sum_{i=1}^{N_{b}}{\gamma_{i}\psi(\lVert \mathbf{x}-\mathbf{s}_{i}\rVert)}$, where $\psi$ is selected as the famous multi-quadric (MQ) RBFs $\psi(r)=\sqrt{r^2+c_{MQ}^2}$ with $c_{MQ}=8$, the RBF centers $\{\mathbf{s}_{i}\} \in \Gamma$ with $N_{b}=4 \sqrt{N}-4$,and the unknown coefficients $\gamma_{i}$ can be determined by satisfying the boundary conditions (\ref{3.27}). Fig. \ref{Fig.2} plots the numerical solutions and the corresponding relative errors of the proposed method with $K=5$ and $N=36$ in Example 8 with $\alpha=1.5$ and $\alpha(t)=1.4+t/10$ in the range $[0,1]$. It can be observed from Fig. \ref{Fig.2} that the largest numerical errors appears in the central region of the computational domain $\Omega$ at $T=1$.

\section{Conclusions}\label{sec con}
This paper presents a novel semi-analytical collocation method to solve multi-term variable-order time fractional partial differential equations (VOTFPDEs). In the proposed numerical scheme it employs the Fourier series expansion for spatial discretization, which transforms the original multi-term VOTFPDEs into a sequence of multi-term variable-order time fractional ordinary differential equations (VOTFODEs). Then these VOTFODEs can be solved by using the recent-developed backward substitution method with the M$\ddot{u}$ntz polynomials. The present numerical experiments verify that such a combination yields an easy-to-implement numerical method that is a competitive alternative for solving the multi-term VOTFPDEs. Numerical investigations show that the proposed BSM with few terms of M$\ddot{u}$ntz polynomials provides satisfactory simulations for the VOTFODEs or the time discretization of VOTFPDEs, and it has high approximation orders with $AO > 4$. In addition, it avoids the error accumulation effect in the standard FDMs. Besides, the proposed numerical results are in good agreement with the analytical solutions and the reference results for solving 1D and 2D VOTFPDEs, and it has rapid convergence rates with $CO>1.8$.

\section*{Acknowledgements}
The work described in this paper was supported by the National Science Funds of China (Grant Nos. 11772119),  the Fundamental Research Funds for the Central Universities (Grant No. B200202124), the State Key Laboratory of Mechanics and Control of Mechanical Structures (Nanjing University of Aeronautics and astronautics) (Grant No. MCMS-E-0519G01), Alexander von Humboldt Research Fellowship (ID: 1195938), the 111 Project (Grant No. B12032) and the Six Talent Peaks Project in Jiangsu Province of China (Grant No. 2019-KTHY-009).



\section*{References}
\bibliographystyle{elsarticle-num}
\bibliography{BSM01}

\begin{thebibliography}{10}
\expandafter\ifx\csname url\endcsname\relax
  \def\url#1{\texttt{#1}}\fi
\expandafter\ifx\csname urlprefix\endcsname\relax\def\urlprefix{URL }\fi
\expandafter\ifx\csname href\endcsname\relax
  \def\href#1#2{#2} \def\path#1{#1}\fi

\bibitem{hilfer2000applications}
R.~Hilfer, Applications of fractional calculus in physics, World Scientific,
  2000.

\bibitem{mandelbrot1968fractional}
B.~B. Mandelbrot, J.~W. Van~Ness, Fractional brownian motions, fractional
  noises and applications, SIAM review 10~(4) (1968) 422--437.

\bibitem{meerschaert2004finite}
M.~M. Meerschaert, C.~Tadjeran, Finite difference approximations for fractional
  advection--dispersion flow equations, Journal of Computational and Applied
  Mathematics 172~(1) (2004) 65--77.

\bibitem{scalas2000fractional}
E.~Scalas, R.~Gorenflo, F.~Mainardi, Fractional calculus and continuous-time
  finance, Physica A: Statistical Mechanics and its Applications 284~(1) (2000)
  376--384.

\bibitem{chechkin2005fractional}
A.~Chechkin, R.~Gorenflo, I.~Sokolov, Fractional diffusion in inhomogeneous
  media, Journal of Physics A: Mathematical and General 38~(42) (2005) L679.

\bibitem{santamaria2006anomalous}
F.~Santamaria, S.~Wils, E.~De~Schutter, G.~J. Augustine, Anomalous diffusion in
  purkinje cell dendrites caused by spines, Neuron 52~(4) (2006) 635--648.

\bibitem{coimbra2003mechanics}
C.~F. Coimbra, Mechanics with variable-order differential operators, Annalen
  der Physik 12~(11-12) (2003) 692--703.

\bibitem{zhao2015second}
X.~Zhao, Z.-z. Sun, G.~E. Karniadakis, Second-order approximations for variable
  order fractional derivatives: algorithms and applications, Journal of
  Computational Physics 293 (2015) 184--200.

\bibitem{zhuang2009numerical}
P.~Zhuang, F.~Liu, V.~Anh, I.~Turner, Numerical methods for the variable-order
  fractional advection-diffusion equation with a nonlinear source term, SIAM
  Journal on Numerical Analysis 47~(3) (2009) 1760--1781.

\bibitem{sun2012finite}
H.~Sun, W.~Chen, C.~Li, Y.~Chen, Finite difference schemes for variable-order
  time fractional diffusion equation, International Journal of Bifurcation and
  Chaos 22~(04) (2012) 1250085.

\bibitem{hashemi2016numerical}
M.~S. Hashemi, D.~Baleanu, Numerical approximation of higher-order
  time-fractional telegraph equation by using a combination of a geometric
  approach and method of line, Journal of Computational Physics 316 (2016)
  10--20.

\bibitem{zhang2012compact}
Y.-n. Zhang, Z.-z. Sun, X.~Zhao, Compact alternating direction implicit scheme
  for the two-dimensional fractional diffusion-wave equation, SIAM Journal on
  Numerical Analysis 50~(3) (2012) 1535--1555.

\bibitem{liu2013numerical}
F.~Liu, M.~Meerschaert, R.~McGough, P.~Zhuang, Q.~Liu, Numerical methods for
  solving the multi-term time-fractional wave-diffusion equation, Fractional
  Calculus and Applied Analysis 16~(1) (2013) 9--25.

\bibitem{langlands2005accuracy}
T.~Langlands, B.~I. Henry, The accuracy and stability of an implicit solution
  method for the fractional diffusion equation, Journal of Computational
  Physics 205~(2) (2005) 719--736.

\bibitem{cui2009compact}
M.~Cui, Compact finite difference method for the fractional diffusion equation,
  Journal of Computational Physics 228~(20) (2009) 7792--7804.

\bibitem{lin2009stability}
R.~Lin, F.~Liu, V.~Anh, I.~Turner, Stability and convergence of a new explicit
  finite-difference approximation for the variable-order nonlinear fractional
  diffusion equation, Applied Mathematics and Computation 212~(2) (2009)
  435--445.

\bibitem{chen2010numerical}
C.-M. Chen, F.~Liu, V.~Anh, I.~Turner, Numerical schemes with high spatial
  accuracy for a variable-order anomalous subdiffusion equation, SIAM Journal
  on Scientific Computing 32~(4) (2010) 1740--1760.

\bibitem{lin2007finite}
Y.~Lin, C.~Xu, Finite difference/spectral approximations for the
  time-fractional diffusion equation, Journal of Computational Physics 225~(2)
  (2007) 1533--1552.

\bibitem{li2011numerical}
C.~Li, Z.~Zhao, Y.~Chen, Numerical approximation of nonlinear fractional
  differential equations with subdiffusion and superdiffusion, Computers \&
  Mathematics with Applications 62~(3) (2011) 855--875.

\bibitem{katsikadelis2011bem}
J.~T. Katsikadelis, The bem for numerical solution of partial fractional
  differential equations, Computers \& Mathematics with Applications 62~(3)
  (2011) 891--901.

\bibitem{brunner2010numerical}
H.~Brunner, L.~Ling, M.~Yamamoto, Numerical simulations of 2d fractional
  subdiffusion problems, Journal of Computational Physics 229~(18) (2010)
  6613--6622.

\bibitem{chen2010fractional}
W.~Chen, L.~Ye, H.~Sun, Fractional diffusion equations by the kansa method,
  Computers \& Mathematics with Applications 59~(5) (2010) 1614--1620.

\bibitem{fu2013boundary}
Z.-J. Fu, W.~Chen, H.-T. Yang, Boundary particle method for laplace transformed
  time fractional diffusion equations, Journal of Computational Physics 235
  (2013) 52--66.

\bibitem{fu2015method}
Z.-J. Fu, W.~Chen, L.~Ling, Method of approximate particular solutions for
  constant-and variable-order fractional diffusion models, Engineering Analysis
  with Boundary Elements 57 (2015) 37--46.

\bibitem{borwein1994muntz}
P.~Borwein, T.~Erd{\'e}lyi, J.~Zhang, M{\"u}ntz systems and orthogonal
  m{\"u}ntz-legendre polynomials, Transactions of the American Mathematical
  Society 342~(2) (1994) 523--542.

\bibitem{esmaeili2011numerical}
S.~Esmaeili, M.~Shamsi, Y.~Luchko, Numerical solution of fractional
  differential equations with a collocation method based on m{\"u}ntz
  polynomials, Computers \& Mathematics with Applications 62~(3) (2011)
  918--929.

\bibitem{mokhtary2016muntz}
P.~Mokhtary, F.~Ghoreishi, H.~Srivastava, The m{\"u}ntz-legendre tau method for
  fractional differential equations, Applied Mathematical Modelling 40~(2)
  (2016) 671--684.

\bibitem{reutskiy2017new}
S.~Y. Reutskiy, A new semi-analytical collocation method for solving multi-term
  fractional partial differential equations with time variable coefficients,
  Applied Mathematical Modelling 45 (2017) 238--254.

\bibitem{chen2014recent}
W.~Chen, Z.~Fu, C.~Chen, Recent advances in radial basis function collocation
  methods, Springer, 2014.

\bibitem{yang2012computationally}
Q.~Yang, T.~J. Moroney, F.~Liu, I.~Turner, Computationally efficient methods
  for solving time-variable-order time-space fractional reaction-diffusion
  equation, In 5th IFAC Symposium on Fractional Differentiation and its
  Applications, 14-17 May 2012, Hohai University, Nanjing, China.

\bibitem{wei2012analysis}
L.~Wei, Y.~He, X.~Zhang, S.~Wang, Analysis of an implicit fully discrete local
  discontinuous galerkin method for the time-fractional schr{\"o}dinger
  equation, Finite Elements in Analysis and Design 59 (2012) 28--34.

\end{thebibliography}


%
%
%
\end{document}